\documentclass[12pt,apaper]{amsart}
\usepackage{amsmath, amsthm, amssymb, enumerate, mathtools}

\usepackage[breaklinks=true,colorlinks=true,linkcolor=blue,citecolor=red,urlcolor=green]{hyperref}

\usepackage[margin=1in,footskip=0.25in]{geometry}

\usepackage{color}

\usepackage[utf8]{inputenc}
\usepackage{eepic,epic}
\usepackage{esint}
\usepackage{multirow}
\usepackage{array}
\usepackage[utf8]{inputenc}
\usepackage{graphicx}
\usepackage{caption}
\usepackage{booktabs}
\usepackage{siunitx}
\usepackage{makecell}
\usepackage{mathrsfs}
\usepackage{enumitem}

\usepackage{theoremref}

\def\to{\longrightarrow}
\def\mapsto{\longmapsto}

\def\cwedge{\bigcirc\kern-1.07em\wedge\ }

\font\medbf=ptmbo at 12pt

\theoremstyle{plain}

\newtheorem{thm}{Theorem}[section]
\newtheorem{lem}[thm]{Lemma}

\newtheorem{cor}[thm]{Corollary}

\newcommand{\bal}{\begin{aligned}}
\newcommand{\eal}{\end{aligned}}

\usepackage{float}

\theoremstyle{definition}
\newtheorem*{defn*}{Definition}

\theoremstyle{remark}
\newtheorem{rem}[thm]{Remark}
\newtheorem*{pf}{Proof}

\numberwithin{equation}{section}

\newcommand{\w}{^{\phantom i}}

\newcommand{\hyp}{\hskip.5pt\vbox
{\hbox{\vrule width2.5ptheight0.5ptdepth0pt}\vskip2pt}\hskip.5pt}
\newcommand{\hs}{\hskip.7pt}
\newcommand{\hh}{\hskip.4pt}
\newcommand{\hn}{\hskip-.4pt}
\newcommand{\nh}{\hskip-.7pt}
\newcommand{\nnh}{\hskip-1pt}
\newcommand{\syt}{\mathrm{a}}
\newcommand{\stt}{\mathrm{b}}
\newcommand{\sym}{\mathrm{b}}
\newcommand{\ric}{\mathrm{r}}

\newcommand{\sca}{\mathrm{s}}
\newcommand{\ein}{\mathrm{e}}
\newcommand{\met}{\mathrm{g}}
\newcommand{\trc}{\mathrm{trc}}

\newcommand{\yq}{\xi}
\def\bbZ{\mathsf{Z\hskip-4ptZ}}
\newcommand{\bbH}{\mathrm{I\!H}}
\newcommand{\bbR}{\mathrm{I\!R}}
\newcommand{\bbC}{{\mathchoice {\setbox0=\hbox{$\displaystyle\mathrm{C}$}
\hbox{\hbox to0pt{\kern0.4\wd0\vrule height0.9\ht0\hss}\box0}} 
{\setbox0=\hbox{$\textstyle\mathrm{C}$}\hbox{\hbox 
to0pt{\kern0.4\wd0\vrule height0.9\ht0\hss}\box0}} 
{\setbox0=\hbox{$\scriptstyle\mathrm{C}$}\hbox{\hbox 
to0pt{\kern0.4\wd0\vrule height0.9\ht0\hss}\box0}} 
{\setbox0=\hbox{$\scriptscriptstyle\mathrm{C}$}\hbox{\hbox 
to0pt{\kern0.4\wd0\vrule height0.9\ht0\hss}\box0}}}}

\begin{document}	
\title[Weakly Einstein curvature tensors] 
	{Weakly Einstein curvature tensors}
%Wooseok Shin
\author[A.\,Derdzinski, J.\,H.\,Park, W\nnh.\,Shin]{Andrzej 
Derdzinski$^{1}$, JeongHyeong Park$^{2}$, Wooseok Shin$^{2}$}
\address{$^{1}$Department of Mathematics, The Ohio State University,
Columbus, OH 43210, USA}
\address{$^{2}$Department of Mathematics, Sungkyunkwan University, Suwon,
16419, Korea}

\email{andrzej@math.ohio-state.edu, parkj@skku.edu, tlsdntjr@skku.edu}

%\thanks {%The research of
%YE  was supported by Basic Science Research Program through the National Research Foundation of Korea(NRF) funded by the Ministry of Education(Grant number RS-2023-00244736.
%The research of
%SK was supported by Basic Science Research Program through the National Research Foundation of Korea(NRF) funded by the Ministry of Education(Grant number RS-2023-00247409.		
%The research of
%NP was partially supported by ARC Discovery grant DP210100951.
%The research of
%JP was supported by the National Research Foundation of Korea(NRF) grant funded by the Korea government(MSIT) (RS-2024-00334956). 	}

\begin{abstract}We classify weakly Einstein algebraic curvature tensors in an 
oriented Euclidean four-space, defined by requiring 
that the three-index contraction of the curvature tensor against itself be 
a multiple of the inner product. This algebraic formulation parallels the
geometric notion of weakly Einstein Riemannian four-manifolds, which include 
con\-for\-mal\-ly flat scalar-flat, and Einstein manifolds. 
Our main result provides a complete classification of non-Einstein weakly
Einstein curvature tensors in dimension four, naturally dividing them into
three disjoint
five-dimensional families of algebraic types. 
These types are explicitly constructed using bases that simultaneously
diagonalize both the Einstein tensor and the (anti)self-dual Weyl tensors,
which consequently proves that such simultaneous diagonalizability follows
from the weakly Einstein property. We also point out that our classification
has immediate applications, and describe how some known 
geometric examples that are neither Einstein, nor con\-for\-mal\-ly flat
scalar-flat (namely, the EPS space and certain K\"ahler surfaces) fit
within our classification framework.
\end{abstract}
	
\maketitle

\section{Introduction}\label{in}
%Following Euh, Park and Sekigawa \cite{EPS13},
One says that a Riemannian four-man\-i\-fold is {\it weakly Ein\-stein\/}
\cite{EPS13} when the three-in\-dex contraction of 
its curvature tensor against itself equals a function times the metric. This
topic has been studied by several authors. Our Refs.\,\cite{AK}, \cite{BV},
\cite{chen}, and \cite{derdzinski-euh-kim-park} -- \cite{marino-villar} list
articles on weakly Ein\-stein manifolds that do not exclude dimension four or
pos\-i\-tive-def\-i\-nite metrics.

The present paper deals with algebraic curvature 
tensors $\,R\,$ in an oriented Euclidean four-space $\,\mathcal{T}\hs$ which
are {\it weakly Ein\-stein\/} in the sense that, analogously, the
three-in\-dex contraction of $\,R\,$ is a multiple of the inner product
$\,\met$. As we point out in Sect.\,\ref{pr},
\begin{equation}\label{iff}
R\mathrm{\ \ is\ weakly\ Ein\-stein\ if\ and\ only\ if\ }\,\,6\hh W\nnh\nh\ein
=-\sca\hh\ein\hh,
\end{equation}
$W\hskip-3pt,\sca\,$ and $\,\ein\,$ denoting, here and below, the Weyl
component
of $\,R$, its scalar curvature, and its Ein\-stein (trace\-less
Ric\-ci) tensor, with $\,W\nnh\nh\ein\,$ standing for the usual action of
algebraic curvature tensors on symmetric $\,(0,2)\,$ tensors. 

Due to (\ref{iff}), the weakly Ein\-stein property 
follows in the two ``trivial'' cases 
\begin{equation}\label{trv}
\ein\,=\,0\hh\mathrm{,\ \ or\ \ }\,W\hn=\,0\,\mathrm{\ and \ }\,\sca\,=\,0\hh,
\end{equation} 
which reflects two known facts about Riemannian four-man\-i\-folds:
they are weakly Ein\-stein whenever they are Ein\-stein
\cite[Sect.\,1.133]{besse}, and, under the assumption of con\-for\-mal
flatness, being weakly Ein\-stein amounts to 
having zero scalar curvature \cite[Theorem 2(i)]{GHMV}; cf.\ also
\cite[the three lines following formula (4.6)]{derdzinski-euh-kim-park}.

Our use of the word ``trivial'' only refers to the
purely algebraic focus of
our discussion: there is nothing trivial in the problem of classifying
Ein\-stein metrics in dimension four.

We provide, in Theorems~\ref{first} --~\ref{third}, an algebraic
classification of nontrivial weakly Ein\-stein curvature tensors, based on
the spectra of the Ein\-stein and Weyl tensors. Our classification leads to 
immediate applications. Namely, the spectral conditions are often easier
to impose, or verify, than the original form of the weakly Ein\-stein
property. As a result, our paper \cite{derdzinski-park-shin} uses
Theorems~\ref{first} --~\ref{third}, rephrased as
\cite[Theorem 5.1 and Remark 5.2]{derdzinski-park-shin}, which includes
\cite[formula (5.3)]{derdzinski-park-shin}, both to show that
certain metrics are weakly Ein\-stein, 
and to classify some special types of weakly Ein\-stein metrics. For the
former, see \cite[Sect.\,10--11]{derdzinski-park-shin}, which use 
\cite[formulae (5.3-b), (5.3-c) and Remark 5.2]{derdzinski-park-shin}; 
for the latter, \cite[Sect.\,19--20]{derdzinski-park-shin}.

Let $\,\Lambda\nnh^\pm$ be the spaces of self-dual and anti-self-dual 
bi\-vec\-tors in our oriented Euclidean four-space $\,\mathcal{T}\nnh$. Any
positive orthonormal basis $\,u_1\w,\dots,u_4\w$ of 
$\,\mathcal{T}\hs$ leads to the length $\,\sqrt{2\,}$ orthogonal
bases of $\,\Lambda\nnh^\pm$ given by
\begin{equation}\label{bas}
u_1\w\wedge u_2\w\pm u_3\w\wedge u_4\w\hh,\quad 
u_1\w\wedge u_3\w\pm u_4\w\wedge u_2\w\hh,\quad 
u_1\w\wedge u_4\w\pm u_2\w\wedge u_3\w\hh.
\end{equation}
The restrictions $\,W^\pm\nnh\nh:\Lambda\nnh^\pm\to\Lambda\nnh^\pm$ of
$\,W\nh$ acting on bi\-vec\-tors are always di\-ag\-o\-nal\-iz\-ed by
some bases of the form (\ref{bas}). See Lemma~\ref{eqvrt}(a). However, one
cannot, in general, expect $\,u_1\w,\dots,u_4\w$ to be in any way related to
the Ein\-stein tensor $\,\ein$.
 
Our first result shows that
the weakly Ein\-stein case is quite exceptional in this regard.
\begin{thm}\label{tzero}Let\/ $\,R\,$ be a weakly Ein\-stein
algebraic curvature tensor in an oriented Euclidean\/ four-space\/
$\,\mathcal{T}\nnh$. Then there exists a positive orthonormal basis\/
$\,u_1\w,\dots,u_4\w$ of\/ $\,\mathcal{T}\nnh$, consisting of eigen\-vec\-tors 
of\/ $\,\ein$, and such that the bases\/ {\rm(\ref{bas})} of\/
$\,\Lambda\nnh^\pm\nnh\hn$ di\-ag\-o\-nal\-ize\/ $\,W^\pm\nh$.
\end{thm}
We easily derive Theorem~\ref{tzero}, at the end of Sect.\,\ref{td}, from
the next three results, which constitute a complete classification of
(non-Ein\-stein) weakly Ein\-stein 
algebraic curvature tensors in dimension four. The resulting
al\-ge\-bra\-ic-e\-quiv\-a\-lence types form three disjoint classes,
depending, first, on whether $\,\ein\,$ does, or does not, have two
distinct double eigen\-val\-ues. The latter case is further split into two
subcases, with $\,\sca=0\,$ and $\,\sca\ne0$. Each of the three disjoint
classes is five-di\-men\-sion\-al, as we construct it using five free
parameters. 

In all three theorems below, $\,u_1\w,\dots,u_4\w$ is a fixed positive 
orthonormal basis of an oriented Euclidean four-space $\,\mathcal{T}\hs$ and,
invoking the Sing\-er-Thorpe theorem \cite[Sect.\,1.126]{besse}, 
%$\,\Lambda\nnh^\pm$ are the spaces of self-dual and anti-self-dual 
%bi\-vec\-tors.
to define an algebraic curvature tensor $\,R\,$ in
$\,\mathcal{T}\nnh$, we prescribe its scalar curvature $\,\sca$, and
\begin{equation}\label{req}
\begin{array}{l}
\mathrm{require\ the\ bases\ }\hs u_1\w,\dots,u_4\w\hn\mathrm{\ of\ 
}\hs\mathcal{T}\nh\mathrm{\ and\ (\ref{bas})\ of\ }\,\Lambda\nnh^\pm\mathrm{\
to\ di\-ag\-o\-nal\-ize\ the\ Ein\-stein}\\
\mathrm{tensor\ }\,\hs\ein\hs\,\mathrm{\ and,\ respectively,\ 
}\,W^\pm\nnh\nh:\Lambda\nnh^\pm\to\Lambda\nnh^\pm\mathrm{,\ with\ some\
specific\ eigen\-values.}
\end{array}
\end{equation}
\begin{rem}\label{equiv}In each of the next three theorems, to define 
$\,W\nh$ one can obviously use item (b) instead of prescribing 
the spectra of $\,W^\pm\nnh$.
\end{rem}
\begin{thm}\label{first}Given\/ %with unordered systems 
$\,\mu_1\w,\mu_2\w,\mu_3\w,\mu_4\w,c_2\w,c_3\w,c_4\w\in\bbR\,$ 
with\/ $\,\mu_1\w+\mu_2\w+\mu_3\w+\mu_4\w=c_2\w+c_3\w+c_4\w=0$, 
let\/ $\,\mu_1\w,\mu_2\w,\mu_3\w,\mu_4\w$ and
$\,\pm c_2\w,\pm c_3\w,\pm c_4\w$ be the eigen\-values in\/ 
{\rm(\ref{req})}. If we set\/ $\,\sca=0$, then
\begin{enumerate}
\item[{\rm(a)}]the resulting algebraic curvature tensor\/ $\,R\,$ is weakly
Ein\-stein,
\item[{\rm(b)}]the only possibly-nonzero components of\/ $\,W\nnh$ in the
frame\/ $\,u_1\w,\dots,u_4\w$ are those algebraically related to 
$\,W_{\!1234}\w,W_{\!1342}\w,W_{\!1423}\w$, where\/ 
$\,(W_{\!1234}\w,W_{\!1342}\w,W_{\!1423}\w)=(c_2\w,c_3\w,c_4\w)$,
\item[{\rm(c)}]every non-Ein\-stein, weakly Ein\-stein algebraic curvature
tensor in dimension four, having zero scalar curvature, arises as above,
unless\/ $\,\ein\,$ has two double eigen\-val\-ues.
\end{enumerate}
\end{thm}
\begin{thm}\label{secnd}If\/ 
$\,\sca,\lambda,\mu,c_2\w,c_3\w,c_4\w\in\bbR\,$ 
have\/ $\,c_2\w+c_3\w+c_4\w=0\,$ and\/ $\,\lambda>\mu\ge0$, and 
we choose the scalar curvature to be\/ $\,\sca$, then, for the
eigen\-values in\/ {\rm(\ref{req})} provided by\/ 
$\,(\mu_1\w,\mu_2\w,\mu_3\w,\mu_4\w)=(-\nnh\lambda,-\nh\mu,\mu,\lambda)\,$
and\/ $\,\pm c_2\w-\sca/12,\,\pm c_3\w-\sca/12,\,\pm c_4\w+\sca/6$,
\begin{enumerate}
\item[{\rm(a)}]the resulting algebraic curvature tensor\/ $\,R\,$ is weakly
Ein\-stein,
\item[{\rm(b)}]the only possibly-nonzero components of\/ $\,W\nnh$ in the
frame\/ $\,u_1\w,\dots,u_4\w$ are those algebraically related to the ones
occurring in the ordered triples\/
$\,(W_{\!1212}\w,W_{\!1313}\w,W_{\!1414}\w)
=(W_{\!3434}\w,W_{\!2424}\w,W_{\!2323}\w)=(-\sca/12,-\sca/12,\sca/6)\,$ and\/
$\,(W_{\!1234}\w,W_{\!1342}\w,W_{\!1423}\w)=(c_2\w,c_3\w,c_4\w)$,
\item[{\rm(c)}]every non-Ein\-stein, weakly Ein\-stein algebraic curvature
tensor in dimension four, with\/ $\,\sca\ne0$, arises as 
above, except when\/ $\,\ein\,$ has two double eigen\-val\-ues.
\end{enumerate}
\end{thm}
\begin{thm}\label{third}For\/ 
$\,\sca,\lambda,\yq,c_2\w,c_3\w,c_4\w\in\bbR\,$ with\/
$\,c_2\w+c_3\w+c_4\w=0\,$ and\/ $\,\lambda>0$, let the scalar curvature be\/
$\,\sca$. With the eigen\-values in\/ {\rm(\ref{req})} given by\/ 
$\,(\mu_1\w,\mu_2\w,\mu_3\w,\mu_4\w)=(-\nnh\lambda,-\nnh\lambda,\lambda,\lambda)$
and\/ $\,\pm c_2\w-\sca/12,\,\pm c_3\w+\yq-\sca/12,\,\pm c_4\w-\yq+\sca/6$,
\begin{enumerate}
\item[{\rm(a)}]the resulting algebraic curvature tensor\/ $\,R\,$ is weakly
Ein\-stein,
\item[{\rm(b)}]the only possibly-nonzero components of\/ $\,W\nnh$ in the
frame\/ $\,u_1\w,\dots,u_4\w$ are algebraically related to those 
in the ordered triples\/
$\,(W_{\!1212}\w,W_{\!1313}\w,W_{\!1414}\w)
=(W_{\!3434}\w,W_{\!2424}\w,W_{\!2323}\w)$ $=(-\sca/12,\yq-\sca/12,-\yq+\sca/6)\,$
and\/ $\,(W_{\!1234}\w,W_{\!1342}\w,W_{\!1423}\w)=(c_2\w,c_3\w,c_4\w)$,
\item[{\rm(c)}]every weakly Ein\-stein algebraic curvature tensor in dimension
four such that\/ $\,\ein$ has two distinct double eigen\-val\-ues 
arises as described above.
\end{enumerate}
\end{thm}
The multiplicities of the eigen\-val\-ues of $\,\ein$, listed in weakly
descending order, form one of the five strings, of which the last four
(the non-Ein\-stein cases) are of interest to us:
\begin{equation}\label{str}
\begin{array}{llllll}
\text{\medbf4}\hh,&\text{\medbf31}\hh,&\text{\medbf22}\hh,
&\text{\medbf211}\hh,&\text{\medbf1111}\hh,\quad&\mathrm{\ so\
that}\\
d=1\hh,\quad&d=2\hh,\quad&d=2\hh,\quad&d=3\hh,\quad&d=4\hh,&
\end{array}
\end{equation}
$d\,$ being the number of distinct eigen\-val\-ues. As the statements of the
last three results clearly indicate, Theorem~\ref{first} involves
the multiplicities $\,\text{\medbf31}$, $\,\text{\medbf211}\,$ and
$\,\text{\medbf1111}$, Theorem~\ref{secnd}  -- only $\,\text{\medbf211}\,$ 
and $\,\text{\medbf1111}$, while Theorem~\ref{third} deals solely with
the multiplicities $\,\text{\medbf22}$.

Our proofs of Theorems~\ref{first} --~\ref{third} are split %divided into
%parts contained in
between Sect.\,\ref{pp},~\ref{ft},~\ref{te} for the first two theorems,
and Sect.\,\ref{pp},~\ref{td} for the last one, as we establish each
item (c) %their final clauses
by cases, based on (\ref{str}).

An algebraic curvature tensor $\,R\,$ in an oriented Euclidean four-space
$\,\mathcal{T}\hs$ is said to be of {\it K\"ah\-ler type\/} if $\,\ein\,$
and $\,W^+\nnh\nh:\Lambda\nnh^+\to\Lambda\nnh^+$ are simultaneously
di\-ag\-o\-nal\-ized by some positive orthonormal basis 
$\,u_1\w,\dots,u_4\w$ of $\,\mathcal{T}\hs$ and the basis (\ref{bas}) for the
$\,+\,$ sign with the respective eigen\-val\-ues 
$\,\mu_1\w,\mu_2\w,\mu_3\w,\mu_4\w$ and $\,(\sca/6,-\sca/12,-\sca/12)\,$ such
that $\,\mu_1\w=\mu_2\w$ and $\,\mu_3\w=\mu_4\w$. This is well known
\cite[pp.\,485--486]{hdg00} to be the case for the curvature tensor of every 
K\"ah\-ler surface.

The following fact is an obvious consequence Theorems~\ref{tzero}, 
~\ref{third} and Remark~\ref{equiv}.
\begin{cor}\label{ktype}For\/ $\,\sca,\lambda,\yq\in\bbR$, 
letting the scalar curvature and %to be\/ $\,\sca$, 
the spectra in\/ {\rm(\ref{req})} be
\begin{equation}\label{spc}
\sca\hh,\quad(-\nnh\lambda,-\nnh\lambda,\lambda,\lambda)\hh,\quad
(\sca/6,-\sca/12,-\sca/12)\hh,\quad
(-\sca/3,2\yq-\sca/12,-\nh2\yq+5\hh\sca/12)\hh,
\end{equation}
%and the scalar curvature to be\/ $\,\sca$, 
we define a K\"ah\-ler-type weakly Ein\-stein algebraic curvature tensor.

Conversely, every non-Ein\-stein K\"ah\-ler-type weakly Ein\-stein algebraic
curvature tensor in dimension four is obtained as described here from some\/
$\,\sca,\yq\,$ and\/ $\,\lambda>0$.

Rather than prescribing the spectra of\/ $\,W^\pm\nnh\nh$ in\/ {\rm(\ref{spc})}
one can equivalently require that the condition\/ {\rm(b)} in
Theorem\/~{\rm\ref{third}} be 
satisfied with\/ $\,(c_2\w,c_3\w,c_4\w)=(\sca/4,-\yq,\yq-\sca/4)$. 
\end{cor}
It is convenient to use the term {\it proper\/} when referring 
to weakly Ein\-stein manifolds or metrics that are neither Ein\-stein,
nor con\-for\-mal\-ly flat and sca\-lar-flat. 
There are three narrow classes of known examples
\cite{derdzinski-euh-kim-park,derdzinski-park-shin,EPS14} of proper weakly
Ein\-stein Riemannian four-man\-i\-folds. 
They have abundant local self-i\-som\-e\-tries: 
some are locally homogeneous \cite{EPS14}, the others have local
co\-ho\-mo\-ge\-ne\-i\-ty one \cite[Remark 12.1]{derdzinski-euh-kim-park} 
or two \cite[Sect. 21]{derdzinski-park-shin}. It turns out that
\begin{equation}\label{rea}
\begin{array}{l}
\mathrm{the\ weakly\ Ein\-stein\ K}\ddot{\mathrm a}\mathrm{h\-ler\ surfaces\
constructed\ in\
}\text{\rm\cite[Sect.\,12]{derdzinski-euh-kim-park}}\mathrm{\ real}\hyp\\
\mathrm{ize\ all\ parameter\ triples\ }\,\hs(\sca,\lambda,\xi)\hs\,\mathrm{\
of\ Corollary\ \ref{ktype}\ for\ which\ }\,\,\sca\hs=\hs8\hs\xi,
\end{array}
\end{equation}
as we show in Sect.\,\ref{kg}, where we explain how the curvature tensors of
the examples in \cite{derdzinski-euh-kim-park,EPS14} fit within our algebraic 
classification. The analogous %The same
question about the examples in \cite{derdzinski-park-shin} 
is answered in \cite[Sect.\,21]{derdzinski-park-shin}.

In \cite[p.\,32]{shapiro} Shapiro establishes a classification theorem for
weakly Ein\-stein algebraic curvature tensors $\,R\,$ (in all dimensions),
where $\,R\,$ is also assumed \cite[Definition 1.3]{shapiro} to be the
Kul\-kar\-ni-No\-mi\-zu square of some symmetric $\,2$-ten\-sor $\,\sym$. 
%$\,R=a\wedge a\,$ for 
Her result states that $\,R\,$ is weakly Ein\-stein if and only if $\,\sym\,$
is either a multiple of the inner product $\,\met$, or has exactly two,
mutually opposite eigen\-values relative to $\,\met$. In dimension four, all
such $\,R$ are ``trivial'' in our sense: $\,d\in\{1,2\}\,$ in (\ref{str}), for
$\,\sym\,$ rather than $\,\ein$, and the possible multiplicities are 
$\,\text{\medbf4}$, $\,\text{\medbf22}$ and $\,\text{\medbf31}$. They all
yield (\ref{trv}) -- the first two have $\,\ein=0$, the third one
$\,W\nnh\nh=0\,$ and $\,\sca=0$. The last sentence is justified in
Remark~\ref{knsqr}.

\section{Preliminaries}\label{pr}
In a Euclidean $\,n$-space $\,\mathcal{T}\hs$ with the inner product $\,\met$, 
any algebraic curvature tensor $\hs R\hs$ acts on arbitrary $\,(0,2)\,$
tensors $\,\stt\,$ by $\,[R\hh\stt]_{ij}\w=R_{ipjq}\w b\hh^{pq}\nh$. This
action preserves (skew)sym\-metry of $\,\stt\,$ and, in view of the Bianchi
identity, becomes $\,2[R\stt]_{ij}\w=R_{ijpq}\w b\hh^{pq}$ when $\,\stt$ 
is skew-sym\-met\-ric \cite[Sect.\,1.114,\,1.131]{besse}. The {\it triple
contraction\/} of $\,R\,$ is the symmetric $\,(0,2)\,$ tensor $\,\trc\hs R\,$
with $\,[\trc\hs R]_{ij}\w\,=\,\,R_{ikpq}\w R_j\w{}^{kpq}\nh$. Our sign
convention about the curvature tensor is such that the Ric\-ci tensor 
has the components $\,R_{ij}\w=R_{ikj}\w{}^k\nh$. Here (but
nowhere else in the paper) we use $\,g$-in\-dex raising and summation over
repeated indices.

As pointed out in \cite[formula (4.6)]{derdzinski-euh-kim-park}, when $\,n=4$,
\begin{equation}\label{trm}
\trc\hs R-2\hh W\nnh\ric-\sca\hh\ric/3\,\mathrm{\
is\ a\ multiple\ of\ }\,\met\hh,
\end{equation}
$W\hskip-3pt,\sca,\ric\,$ being the Weyl component, scalar curvature, and
Ric\-ci contraction of $\,R$. This trivially implies
(\ref{iff}). In \cite[Theorem 5.1(c)]{derdzinski-euh-kim-park}, (\ref{iff}) is
established for K\"ah\-ler surfaces.

For the remainder of this section, $\,R\,$ is an algebraic curvature 
tensor in an oriented Euclidean four-space $\,\mathcal{T}\nnh$, with
$\,W\hskip-3pt,\sca\,$ as above and the Ein\-stein tensor
$\,\ein=\ric-\sca\hs\met/4$. 

Let a positive orthonormal basis $\,u_1\w,\dots,u_4\w$ of $\,\mathcal{T}\hs$
di\-ag\-o\-nal\-ize $\,\ein\,$ with some eigen\-val\-ues
$\,\mu_1\w,\mu_2\w,\mu_3\w,\mu_4\w$.
Since $\,W\hs$ is well known \cite[p.\,647]{hdg00} to commute with the 
Hodge star, for the components of $\,W\nh$ in the frame $\,u_1\w,\dots,u_4\w$
one has
\begin{equation}\label{wij}
W_{\!ijkl}\w=W_{\!pqrs}\w\mathrm{\ if\ }\,\,(i,j,p,q)\,\mathrm{\ and\
}\,(k,l,r,s)\,\mathrm{\ are\ even\ permutations\ of\ }\,(1,2,3,4).
\end{equation}
Due to (\ref{iff}), 
the weakly Ein\-stein property of $\,R\,$ is equivalent to
\begin{equation}\label{eqv}
\mathrm{a)}\hskip6ptW_{\!ikjk}\w\hs(\mu_k\w-\mu_l\w)=0\hh,\qquad
\mathrm{b)}\hskip6ptW_{\!ijij}\w\hs\mu_j\w+W_{\!ikik}\w\hs\mu_k\w
+W_{\!ilil}\w\hs\mu_l\w=-\sca\hh\mu_i\w/6\hh,
\end{equation}
whenever $\,\{i,j,k,l\}=\{1,2,3,4\}$. Here (\ref{eqv}-a) arises from 
$\,W_{\!ikjk}\w\hs\mu_k\w+W_{\!iljl}\w\hs\mu_l\w=0\,$ and (\ref{wij}): 
$\,W_{\!iljl}\w=-W_{\!ikjk}\w$ since the permutations
$\,(i,l,j,k),(j,l,i,k)\,$ have opposite parities. Let us set
\begin{equation}\label{set}
(a_2\w,a_3\w,a_4\w)=(W_{\!1212}\w,W_{\!1313}\w,W_{\!1414}\w)
=(W_{\!3434}\w,W_{\!2424}\w,W_{\!2323}\w)\hh,
\end{equation}
the second equality being obvious from (\ref{wij}). With (\ref{set}),
equation (\ref{eqv}-b) reads
\begin{equation}\label{mtr}
\left[\begin{matrix}
\mu_2\w&\mu_3\w&\mu_4\w\cr
\mu_1\w&\mu_4\w&\mu_3\w\cr
\mu_4\w&\mu_1\w&\mu_2\w\cr
\mu_3\w&\mu_2\w&\mu_1\w\end{matrix}\right]
\left[\begin{matrix}
a_2\w\cr
a_3\w\cr
a_4\w\end{matrix}\right]
=-\frac{\sca}6\left[\begin{matrix}
\mu_1\w\cr
\mu_2\w\cr
\mu_3\w\cr
\mu_4\w\end{matrix}\right]\nnh.
\end{equation}
To simplify our discussion, here and in Sect.\,\ref{sa}, we let
\begin{equation}\label{rho}
\rho\in\{0,1,2,3\}\,\mathrm{\ denote\ the\ rank\ of\ the\
}\,4\times3\,\mathrm{\ matrix\ in\ (\ref{mtr}).}
\end{equation}
We will encounter some cases where, requiring that
$\,\mu_1\w\le\mu_2\w\le\mu_3\w\le\mu_4\w$, one has
\begin{equation}\label{lbm}
(\mu_1\w,\mu_2\w,\mu_3\w,\mu_4\w)\,=\,(-\nnh\lambda,-\nh\mu,\mu,\lambda)\,\mathrm{\
 with\ }\,\lambda\ge\mu\ge0\,\mathrm{\ and\ }\,\lambda>0\hh.
\end{equation}
\begin{lem}\label{solns}Under the assumption\/ {\rm(\ref{lbm})}, the only 
solutions\/ $\,(a_2\w,a_3\w,a_4\w)\,$ of\/ {\rm(\ref{mtr})} such that\/ 
$\,a_2\w+a_3\w+a_4\w=0\,$ are
\begin{enumerate}
\item[{\rm(a)}]just\/ $\,(a_2\w,a_3\w,a_4\w)=(-\sca/12,-\sca/12,\sca/6)$,  
when\/ $\,\lambda>\mu$,
\item[{\rm(b)}]$(a_2\w,a_3\w,a_4\w)=(-\sca/12,\yq-\sca/12,-\yq+\sca/6)\,$ for
all\/ $\,\yq\in\bbR$, if\/ $\,\lambda=\mu$. 
\end{enumerate}
\end{lem}
\begin{pf}Both (a)\hs--\hs(b) define $\,(a_2\w,a_3\w,a_4\w)\,$ with 
(\ref{mtr}) and $\,a_2\w+a_3\w+a_4\w=0$. If $\,\lambda>\mu\,$ (or,
$\,\lambda=\mu$), then $\rho=2\,$ (or, $\rho=1$) in (\ref{rho}), and
the kernel of the matrix operator in (\ref{mtr}) is spanned by
$\,(1,1,0)\,$ or, respectively, by $\,(1,1,0)\,$ and $\,(0,1,-\nh1)$, so that
the only vectors $\,(x_2\w,x_3\w,x_4\w)\,$ in the kernel having 
$\,x_2\w+x_3\w+x_4\w=0\,$ are $\,(0,0,0)\,$ or, respectively, multiples of
$\,(0,1,-\nh1)$.
\qed
\end{pf}
\begin{rem}\label{diagz}In an oriented Euclidean plane, any trace\-less 
symmetric $\,(0,2)\,$ ten\-sor $\,\syt\,$ has, in a suitable positive
orthonormal basis, a matrix with zeros on the diagonal: we start with
a positive orthonormal basis $\,v,w\,$ di\-ag\-o\-nal\-izing $\,\syt$, and
replace it with $\,(v-w)/\nh\sqrt{2\,},(v+w)/\nh\sqrt{2\,}$.
\end{rem}
  
\section{Proofs of parts {\rm(a)} and {\rm(b)} in Theorems~\ref{first},
~\ref{secnd} and~\ref{third}}\label{pp}
In Theorem~\ref{first}, $\,W\nnh$ sends
$\,u_1\w\wedge u_2\w,\,u_1\w\wedge u_3\w,\,u_1\w\wedge u_4\w, 
u_3\w\wedge u_4\w,\,u_4\w\wedge u_2\w,\,u_2\w\wedge u_3\w$ to, respectively, 
$\,c_2\w u_3\w\wedge u_4\w,\,c_3\w u_4\w\wedge u_2\w,\,c_4\w u_2\w\wedge u_3\w,
c_2\w u_1\w\wedge u_2\w,\,c_3\w u_1\w\wedge u_3\w,\,c_4\w u_1\w\wedge u_4\w$,
and (b) trivially follows. As $\,\sca=0$, this yields (\ref{eqv}), and hence
(a).

Theorem~\ref{secnd}: $\,W\nnh$ sends
$\,u_1\w\wedge u_2\w,\,u_1\w\wedge u_3\w,\,u_1\w\wedge u_4\w, 
u_3\w\wedge u_4\w,\,u_4\w\wedge u_2\w,\,u_2\w\wedge u_3\w$ to, respectively,
\[
\begin{array}{l}
c_2\w u_3\w\wedge u_4\w-\sca\hh u_1\w\wedge u_2\w/12,\,
c_3\w u_4\w\wedge u_2\w-\sca\hh u_1\w\wedge u_3\w/12,\,
c_4\w u_2\w\wedge u_3\w+\sca\hh u_1\w\wedge u_4\w/6\hh,\\
\hskip35ptc_2\w u_1\w\wedge u_2\w-\sca\hh u_3\w\wedge u_4\w/12,\,
c_3\w u_1\w\wedge u_3\w-\sca\hh u_4\w\wedge u_2\w/12,\,
c_4\w u_1\w\wedge u_4\w+\sca\hh u_2\w\wedge u_3\w/6\hh,
\end{array}
\]
proving (b); (a) follows from (\ref{set}), as 
(\ref{eqv}-a) is obvious and Lemma~\ref{solns}(a) yields (\ref{eqv}-b).

Theorem~\ref{third}: $\,W\nnh$ sends
$\,u_1\w\wedge u_2\w,\,u_1\w\wedge u_3\w,\,u_1\w\wedge u_4\w, 
u_3\w\wedge u_4\w,\,u_4\w\wedge u_2\w,\,u_2\w\wedge u_3\w$ to, respectively,
\[
\begin{array}{l}
c_2\w u_3\w\wedge u_4\w-\sca\hh u_1\w\wedge u_2\w/12,\,
c_3\w u_4\w\wedge u_2\w+(\yq-\sca/12)\hh u_1\w\wedge u_3\w,\,
c_4\w u_2\w\wedge u_3\w-(\yq-\sca/6)\hh u_1\w\wedge u_4\w\hh,\\
\hskip15ptc_2\w u_1\w\wedge u_2\w-\sca\hh u_3\w\wedge u_4\w/12,\,
c_3\w u_1\w\wedge u_3\w+(\yq-\sca/12)\hh u_4\w\wedge u_2\w,\,
c_4\w u_1\w\wedge u_4\w-(\yq-\sca/6)\hh u_2\w\wedge u_3\w\hh,
\end{array}
\]
which proves (b); 
with (\ref{eqv}-a) obvious, Lemma~\ref{solns}(b) 
yields (\ref{eqv}-b) and, by (\ref{set}), (a). 

\section{Some algebraic facts}\label{sa}
As before, $\,R\,$ is an algebraic curvature 
tensor in an oriented Euclidean four-space $\,\mathcal{T}\nnh$. We use the
notation of Sect.\,\ref{pr}. Consider the following three conditions:
\begin{equation}\label{con}
\begin{array}{rl}
\mathrm{(i)}&\sca\,\ne\,0\hh,\\
\mathrm{(ii)}&\rho\,<\hs3\,\mathrm{\ in\ (\ref{rho}),}\\
\mathrm{(iii)}&\mathrm{the\ set\
}\,\{|\mu_1\w|,|\mu_2\w|,|\mu_3\w|,|\mu_4\w|\}\,\mathrm{\ has\ at\ most\
two\ elements.}
\end{array}
\end{equation}
\begin{lem}\label{msttw}The condition\/ {\rm(\ref{con}-i)} implies\/ 
{\rm(\ref{con}-ii)}, and\/ {\rm(\ref{con}-ii)} implies\/ 
{\rm(\ref{con}-iii)}. Also,
\begin{enumerate}
\item[{\rm(a)}]if the set\/ $\,\{|\mu_1\w|,|\mu_2\w|,|\mu_3\w|,|\mu_4\w|\}\,$
has more than two elements, then\/ $\,\rho=3$,
\item[{\rm(b)}]if the multiplicities  in\/ {\rm(\ref{str})} are\/ 
$\,\text{\medbf31}$, then\/ $\,\rho=3$,  
\item[{\rm(c)}]if\/ $\,\rho=3$, then\/ $\,\sca=a_2\w=a_3\w=a_4\w=0\,$ in\/
{\rm(\ref{mtr})}.
\end{enumerate}
\end{lem}
\begin{pf}The formula $\,(a_2\w,a_3\w,a_4\w)=(\sca/6,\sca/6,\sca/6)\,$ defines
a solution of (\ref{mtr}). Since, by (\ref{set}), our
$\,(a_2\w,a_3\w,a_4\w)\,$ has
$\,a_2\w+a_3\w+a_4\w=0$, the assumption that $\,\sca\ne0\,$ leads to
nonuniqueness of solutions, and hence (\ref{con}-ii). Next, 
(\ref{con}-ii) amounts to the determinant conditions
\[
\begin{array}{l}
2(\mu_1^3+\mu_2\w\mu_3\w\mu_4\w)=\mu_1\w(\mu_1^2+\mu_2^2+\mu_3^2+\mu_4^2)\hh,\\
2(\mu_2^3+\mu_1\w\mu_3\w\mu_4\w)=\mu_2\w(\mu_1^2+\mu_2^2+\mu_3^2+\mu_4^2)\hh,\\2(\mu_3^3+\mu_1\w\mu_2\w\mu_4\w)=\mu_3\w(\mu_1^2+\mu_2^2+\mu_3^2+\mu_4^2)\hh,\\2(\mu_4^3+\mu_1\w\mu_2\w\mu_3\w)=\mu_4\w(\mu_1^2+\mu_2^2+\mu_3^2+\mu_4^2)\hh,
\end{array}
\]
where it suffices to verify just one equality, as the other three then follow
from invariance of (\ref{eqv}-b) under permutations of $\{1,2,3,4\}$.
Multiplying the four lines above by, respectively, 
$\,\mu_1\w,\mu_2\w,\mu_3\w,\mu_4\w$,
we see that each $\,\sigma=\mu_i^2$ is a root
of the quadratic equation $\,2\sigma^2\nh-p\sigma+q=0$, for
$\,p=\mu_1^2+\mu_2^2+\mu_3^2+\mu_4^2$ and $\,q=2\mu_1\w\mu_2\w\mu_3\w\mu_4\w$,
so that (\ref{con}-iii) follows.

Now (a) is obvious, for (b) we may require that 
$\,(\mu_1\w,\mu_2\w,\mu_3\w,\mu_4\w)\,=\,(\mu,\mu,\mu,-3\mu)\,$ with some
$\,\mu\ne0$, which gives $\rho=3\,$ in (\ref{rho}), while for (c) we have 
$\,(\rho,\sca)=(3,0)$, and the resulting injectivity of the matrix operator
in (\ref{mtr}) yields $\,(a_2\w,a_3\w,a_4\w)=(0,0,0)$.
\qed
\end{pf}

\section{Step one toward the proofs of Theorems~\ref{first} 
--~\ref{third}}\label{so}
The components of $\,W\hh$ in a positive orthonormal basis
$\,u_1\w,\dots,u_4\w$, di\-ag\-o\-nal\-izing $\,\ein\,$ with some
eigen\-val\-ues $\,\mu_1\w,\mu_2\w,\mu_3\w,\mu_4\w$,
are of three kinds:
\begin{equation}\label{cpt}
\mathrm{a)}\hskip6ptW_{\!ijik}\w\hh,\qquad
\mathrm{b)}\hskip6ptW_{\!ijij}\w\hh,\qquad
\mathrm{c)}\hskip6ptW_{\!ijkl}\w\hh,
\end{equation}
where $\,\{i,j,k,l\}=\{1,2,3,4\}$. According to (\ref{eqv}), the three kinds
are subject to three separate conditions: (\ref{eqv}-a) for (\ref{cpt}-a) and 
(\ref{eqv}-b) for (\ref{cpt}-b), with no conditions at all for (\ref{cpt}-c),
so that {\it the components\/} (\ref{cpt}-c) {\it are completely arbitrary\hs}:
\begin{equation}\label{arb}
(W_{\!1234}\w,W_{\!1342}\w,W_{\!1423}\w)=(c_2\w,c_3\w,c_4\w)\,\mathrm{\ with\
}\,c_2\w+c_3\w+c_4\w=0\hh.
\end{equation}
As for (\ref{cpt}-a), equation (\ref{eqv}-a) states that some of them
vanish, and the others are arbitrary.
\begin{lem}\label{oooot}For the multiplicities\/ 
$\,\text{\medbf1111},\,\,\text{\medbf211}\,$ and\/ $\,\text{\medbf31}\,$ in\/ 
{\rm(\ref{str})}, if\/ $\,\mu_1\w\le\mu_2\w\le\mu_3\w\le\mu_4\w$,
\begin{equation}\label{wot}
(W_{\!1212}\w,W_{\!1313}\w,W_{\!1414}\w)
=(W_{\!3434}\w,W_{\!2424}\w,W_{\!2323}\w)=(-\sca/12,-\sca/12,\sca/6)\hh,
\end{equation}
while in the cases\/ $\,\text{\medbf1111}\,$ and\/
$\,\text{\medbf211}$, when\/ $\,\sca\ne0$, one must have\/ {\rm(\ref{lbm})} 
with\/ $\,\lambda>\mu\ge0$.
\end{lem}
\begin{pf}If the multiplicities are $\,\text{\medbf31}$, we get 
(\ref{wot})
from Lemma~\ref{msttw}(b)\hs--\hs(c) and (\ref{set}). 
When $\,\{|\mu_1\w|,\nh|\mu_2\w|,\nh|\mu_3\w|,\nh|\mu_4\w|\}\,$ has
more than two
elements, $\,(\sca,\hn a_2\w,\hn a_3\w,\hn a_4\w)=(0,\hn0,\hn0,\hn0)\,$ by
Lemma~\ref{msttw}(a),\hs(c), and (\ref{wot}) follows from (\ref{set}). Let
$\,\{|\mu_1\w|,|\mu_2\w|,|\mu_3\w|,|\mu_4\w|\}\,$ now have at most 
two elements. For the multiplicities  
$\,\text{\medbf1111}\,$ or $\,\text{\medbf211}$, this clearly implies 
(\ref{lbm}) with $\,\lambda>\mu>0$ or, respectively, $\,\lambda>\mu=0$,
and, due to (\ref{set}), Lemma~\ref{solns}(a) gives (\ref{wot}).

Finally, if one has $\,\sca\ne0$,
and $\,\text{\medbf1111}\,$ or $\,\text{\medbf211}$, the first line of
Lemma~\ref{msttw} gives (\ref{con}-iii), and the preceding two sentences
yield the last line of the lemma.
\qed
\end{pf}

\section{Step two: four or three different eigen\-val\-ues}\label{ft}
We are assuming the multiplicities in (\ref{str}) to be 
$\,\text{\medbf1111}\,$ or $\,\text{\medbf211}$.

For $\,\text{\medbf1111}$, (\ref{eqv}-a) gives $\,W_{\!ikjk}\w=0\,$
whenever $\,\{i,j,k,l\}=\{1,2,3,4\}$. Thus, if $\,\sca=0$ (or,
$\,\sca\ne0$), (\ref{arb}) -- (\ref{wot}) imply (b) in
Theorem~\ref{first} (or, respectively, in Theorem~\ref{secnd}).

Due to Remark~\ref{equiv} and Sect.\,\ref{pp}, along with
the final clause of Lemma~\ref{oooot}, this proves 
Theorems~\ref{first} --~\ref{secnd} when the multiplicities are
$\,\text{\medbf1111}$.

In the case $\,\text{\medbf211}$, 
for a positive orthonormal basis $\,u_1\w,\dots,u_4\w$ di\-ag\-o\-nal\-izing 
the Ein\-stein tensor $\,\ein\,$ with some eigen\-val\-ues 
$\,\mu_1\w,\mu_2\w,\mu_3\w,\mu_4\w$, let us order
$\,\mu_1\w,\mu_2\w,\mu_3\w,\mu_4\w$ in such a way that
\begin{equation}\label{mtn}
\begin{array}{l}
\mu_3\w\,\ne\,\mu_1\w\,\ne\,\mu_2\w\,\ne\,\mu_3\w\,=\,\mu_4\w\hh\mathrm{,\
and\ so,\ by\ (\ref{wij})\ }-\mathrm{(\ref{eqv}}\hyp\mathrm{a),}\\
W_{\!1323}\w=W_{\!1442}\w=b\,\mathrm{,\ for\ some\  }\,b\in\bbR\hh,\\
W_{\!1213}\w=W_{\!1214}\w=W_{\!1223}\w=W_{\!1224}\w=W_{\!1314}\w\\
\hskip30.5pt=W_{\!1334}\w=W_{\!1434}\w=W_{\!2324}\w=W_{\!2334}\w=W_{\!2434}\w
=\,0\hh.
\end{array}
\end{equation}
Thus, 
$\,W\nh$ has in the basis 
$\,u_1\w\wedge u_2\w,\,u_1\w\wedge u_3\w,\,u_1\w\wedge u_4\w,
u_3\w\wedge u_4\w,\,u_4\w\wedge u_2\w,\,u_2\w\wedge u_3\w$ the matrix
\[
\left[\begin{matrix}
-\sca/12&0&0&c_2\w&0&0\cr
0&-\sca/12&0&0&c_3\w&b\cr
0&0&\sca/6&0&b&c_4\w\cr
c_2\w&0&0&-\sca/12&0&0\cr
0&c_3\w&b&0&-\sca/12&0\cr
0&b&c_4\w&0&0&\sca/6\end{matrix}\right]\nnh\nh,
\]
cf.\ (\ref{arb}), the matrices of
$\,W^\pm\nh$ in the bases (\ref{bas}) of $\,\Lambda\nnh^\pm$ are
\begin{equation}\label{wpm}
\left[\begin{matrix}
\pm c_2\w-\sca/12&0&0\cr
0&\pm c_3\w-\sca/12&\pm b\cr
0&\pm b&\pm c_4\w+\sca/6\cr\end{matrix}\right]\nnh\nh,
\end{equation}
and the matrix of $\,\syt=W(u^1\otimes u^2\nh+u^2\otimes u^1)\,$ in the
basis $\,u_1\w,\dots,u_4\w$ is
\begin{equation}\label{mae}
\left[\begin{matrix}
0&\sca/12&0&0\cr
\sca/12&0&0&0\cr
0&0&2b&c_4\w-c_3\w\cr
0&0&c_4\w-c_3&-2b\w\cr
\end{matrix}\right]\nnh\nh,
\end{equation}
$u^1\nh,\dots,u^4$ being the basis of $\,\mathcal{T}\hh^*$ dual to
$\,u_1\w,\dots,u_4\w$. By suitably rotating $\,u_3\w,u_4\w$, while 
leaving $\,u_1\w,u_2\w$ and $\,u^1\nh,u^2$ unchanged, we may require
(Remark~\ref{diagz}) that, in addition, $\,b=0\,$ in (\ref{mae}). Now 
(\ref{wpm}) and the last line of Lemma~\ref{oooot} establish 
Theorems~\ref{first} --~\ref{secnd} for $\,\text{\medbf211}$.

\section{Self-duality in terms of quaternions}\label{sd}
Since a connected Lie group has no open proper subgroups, a homomorphism
between two connected Lie groups of the same dimension, having a finite 
kernel, must be a covering projection. 
Consequently, 
for the group $\,S^3\nh=\mathrm{SU}(2)\subseteq\bbH\,$ of unit quaternions 
we have the two-fold covering homomorphisms
\begin{equation}\label{tfc}
S^3\nh\to\mathrm{SO}(1^\perp\nh),\qquad
S^3\nnh\times S^3\nh\to\mathrm{SO}(\bbH),\qquad
\mathrm{SO}(\bbH)\to\mathrm{SO}(1^\perp\nh)\times\mathrm{SO}(1^\perp\nh)
\end{equation}
given by $\,p\mapsto[1^\perp\ni x\mapsto p\hs x\overline p\in1^\perp\nh]\,$
and, respectively, 
$\,(p,q)\mapsto[\bbH\ni x\mapsto p\hs x\overline q\in\bbH]$,
$\,\pm(p,q)\mapsto(\pm p,\pm q)$. 
Here $\,1^\perp\approx\bbR^3$ is the three-space of pure quaternions. (In fact,
as the multiplicative center of $\,\bbH\,$ is $\,\bbR$, all three 
homomorphisms have the kernel $\,\bbZ_2\w$.)

The two-to-one homomorphisms (\ref{tfc}) justify the notation,
with $\,p,q\in S^3\nh$,
\begin{equation}\label{not}
\pm\nnh\nh(p,q)\in\mathrm{SO}(\bbH)\hh,\qquad \pm p\in\mathrm{SO}(1^\perp\nh)\hh.
\end{equation}
For an oriented Euclidean four-space $\,\mathcal{T}\hs$ with the inner 
product $\,\met$, we use the identification
\begin{equation}\label{ide}
\mathcal{T}^{\wedge2}\hn=\hs\mathfrak{so}\hs(\mathcal{T})
%=[\mathcal{T}\hh^*]^{\wedge2}
\end{equation}
provided by $\,\met$. Thus, $\,u\wedge v$, for
$\,u,v\in\mathcal{T}\nnh$, becomes the endomorphism of $\,\mathcal{T}\hs$
given by
\begin{equation}\label{uvw}
w\,\mapsto\,(u\wedge v)w\,=\,\met(u,w)v\,-\,\met(v,w)u\hh.
\end{equation}
Both summands $\,\Lambda\nnh^\pm\nnh$
of $\,\mathcal{T}^{\wedge2}$ are canonically oriented by the bases 
(\ref{bas}): as $\,\mathrm{SO}(4)\,$ is connected, these bases form
connected sets.
\begin{lem}\label{eqvrt}Let $\,\mathcal{T}\hh$ be an oriented Euclidean\/
four-space.
\begin{enumerate}
\item[{\rm(a)}]Every pair of length $\,\sqrt{2\,}$ positive orthogonal 
bases of\/ $\,\Lambda\nnh^\pm\nnh$ has the form\/ {\rm(\ref{bas})} for %some
a positive orthonormal basis $\,u_1\w,\dots,u_4\w$ of $\,\mathcal{T}\nnh$,
%which is
unique up to an overall sign change.
\item[{\rm(b)}]The action of $\,\mathrm{SO}\hh(u^\perp\nh)\,$ on bases 
$\,u_1\w,\dots,u_4\w$ in {\rm(a)} having a fixed first vector\/ $\,u_1\w=u$ 
results in the standard $\,\mathrm{SO}\hh(1^\perp\nh)\,$ actions on the two
bases\/ {\rm(\ref{bas})}, for the two signs\/ $\,\pm\hh$.
\item[{\rm(c)}]The $\,\hs\mathrm{SO}\hh(2)\nnh\times\mathrm{SO}\hh(2)\,$
action on a positive orthonormal basis $\,u_1\w,\dots,u_4\w$ of
$\,\mathcal{T}\nnh$, which
independently rotates the pairs $\,u_1\w,u_2\w$ and\/ $\,u_3\w,u_4\w$, leads 
to the $\,\hs\mathrm{SO}\hh(2)\nnh\times\mathrm{SO}\hh(2)\,$ action on\/
{\rm(\ref{bas})}, independently rotating the last two bi\-vectors.
\end{enumerate}
\end{lem}
\begin{pf}We may set $\,\mathcal{T}\nh=\bbH=\bbR^4\nh$. Due to (\ref{tfc}), 
any positive orthonormal basis $\,u_1\w,\dots,u_4\w$ of $\,\mathcal{T}\hs$
equals $\,p\overline q,p\hs i\overline q,pj\overline q,p\hh k\overline q\,$ 
for a unique $\,\pm(p,q)\,$ with $\,p,q\in S^3\nh$, cf.\ (\ref{not}). 
The identification (\ref{uvw}) clearly 
turns each $\,u_a\w\wedge u_b\w\pm u_c\w\wedge u_d\w$ in (\ref{bas}) into 
\begin{equation}\label{acs}
\mathrm{a\ complex\ structure\ in\ }\,\mathcal{T}\nnh\mathrm{,\ sending\
}\,u_a\w\mathrm{\ to\ }\,u_b\w\mathrm{\ and\ }\,u_c\w\mathrm{\ to\
}\,\pm\hs u_d\w.
\end{equation}
As a consequence of (\ref{acs}), with the identification (\ref{ide}),
\begin{equation}\label{lrm}
\begin{array}{l}
\mathrm{the\ plus}\hh\hyp\hh\mathrm{sign\ (or,\ minus}\hh\hyp\hh\mathrm{sign)\
triple\ in\ (\ref{bas})\ consists\ of}\\
\mathrm{the\ left\ quaternion\ multiplications\ by\ 
}\,p\hs i\overline p,pj\overline p,p\hh k\overline p\,\mathrm{\ (or,}\\
\mathrm{re\-spec\-tively,\ the\ right\ multiplications\ by\ 
}\,qi\overline q,qj\overline q,qk\overline q\mathrm{).}
\end{array}
\end{equation}
Thus, for the last two-fold covering homomorphism 
$\,\pm(p,q)\mapsto(\pm p,\pm q)\,$ in (\ref{tfc}),
\begin{equation}\label{eqt}
\begin{array}{l}
\mathrm{the\ mapping\ that\ sends\ }\,u_1\w,\dots,u_4\w\,\mathrm{\ to\
(\ref{bas})\ is}\\
\mathrm{equivar\-i\-ant\ relative\ to\
}\,\pm\nnh\nnh(p,q)\,\mapsto\hs(\pm p,\pm q)\hh.
\end{array}
\end{equation}
This proves (a). 
For (b), we may require that, under our
identification $\,\mathcal{T}\nh=\bbH$, the vector $\,u_1\w$ correspond
to $\,1$. The iso\-tropy group $\,\mathrm{SO}\hh(u^\perp\nh)\,$ in (b) thus
consists of $\,\pm(p,q)\in\mathrm{SO}\hh(\bbH)$ with $\,p=q$, and (b) 
follows from (\ref{eqt}).

Finally, (c) becomes obvious if we note that replacing $\,(p,q)\,$ with
$\,(pe^{i\phi}\nh,qe^{-i\theta})$, where $\,e^{i\phi}=\cos\phi+i\sin\phi$,
subjects the first (or, second) pair in the quadruple 
$\,(p\overline q,p\hs i\overline q,pj\overline q,p\hh k\overline q)\,$ to
a rotation by the angle $\,\phi+\theta\,$ (or, $\,\phi-\theta$) in
the plane spanned by it while, at the same time,
the last two vectors in the triples
$\,p\hs i\overline p,pj\overline p,p\hh k\overline p\,$ and
$\,\,qi\overline q,qj\overline q,qk\overline q\,$ 
undergo rotations by the angles $\,2\phi\,$ and $\,-\nh2\theta$.
\qed
\end{pf}
Under the identifications (\ref{ide}) for $\,\mathcal{T}\nh=\bbH$, 
(\ref{lrm}) clearly implies that
\begin{equation}\label{lpl}
\begin{array}{l}
\Lambda\nnh^+\nnh\nnh\mathrm{\ (or,\ }\Lambda\nnh^-\nh\mathrm{)\ equals\ the\
 space\ of\ all\ left\ (or,\nnh\ right)}\\
\mathrm{quaternion\hs\hs\ multiplications\hs\hs\ by\hs\hs\ pure\hs\hs\ 
quaternions.}
\end{array}
\end{equation}
%\begin{rem}\label{lpmlr}
%\end{rem}

\section{Step three: one triple and one simple eigen\-val\-ue}\label{te}
The multiplicities in (\ref{str}) are now $\,\text{\medbf31}$, and we are free
to assume that
\begin{equation}\label{mon}
\mu_1\w\,\ne\,\mu_2\w\,=\,\mu_3\w\,=\,\mu_4\w\hh.
\end{equation}
Lemma~\ref{msttw}(b)\hs--\hs(c) gives $\,\sca=a_2\w=a_3\w=a_4\w=0\,$ in 
(\ref{mtr}).

By default, this -- combined with Sect.\,\ref{pp} and~\ref{ft}, completes the
proof of Theorem~\ref{secnd}.

According to (\ref{set}), it also clearly yields
\[%\begin{equation}\label{mtw}
(W_{\!1212}\w,W_{\!1313}\w,W_{\!1414}\w,
W_{\!3434}\w,W_{\!2424}\w,W_{\!2323}\w)=(0,0,0,0,0,0)\hh,
\]%\end{equation}
while, by (\ref{wij}) -- (\ref{eqv}-a) and (\ref{mon}), 
for some $\,b_2\w,b_3\w,b_4\w$,
\[%\begin{equation}\label{smb}
\begin{array}{l}
(W_{\!1242}\w,W_{\!1223}\w,W_{\!1323}\w)=
(W_{\!1334}\w,W_{\!1434}\w,W_{\!1442}\w)=(b_4\w,b_3\w,b_2\w)\hh,\\
(W_{\!1213}\w,W_{\!1214}\w,W_{\!1314}\w,
W_{\!2324}\w,W_{\!2334}\w,W_{\!2434}\w)=(0,0,0,0,0,0)\hh.
\end{array}
\]%\end{equation}
Thus, the matrix of %consequently,
$\,W\nh$ in the basis 
$\,u_1\w\wedge u_2\w,\,u_1\w\wedge u_3\w,\,u_1\w\wedge u_4\w,
u_3\w\wedge u_4\w,\,u_4\w\wedge u_2\w,\,u_2\w\wedge u_3\w$ is
\[
\left[\begin{matrix}
0&0&0&c_2\w&b_4\w&b_3\w\cr
0&0&0&b_4\w&c_3\w&b_2\w\cr
0&0&0&b_3\w&b_2\w&c_4\w\cr
c_2\w&b_4\w&b_3\w&0&0&0\cr
b_4\w&c_3\w&b_2\w&0&0&0\cr
b_3\w&b_2\w&c_4\w&0&0&0\end{matrix}\right]\nnh\nh,
\]
cf.\ (\ref{arb}), and the matrices of
$\,W^\pm\nh$ in the bases (\ref{bas}) of $\,\Lambda\nnh^\pm$ are
\[%\begin{equation}\label{mwp}
\pm\nnh\left[\begin{matrix}
c_2\w&b_4\w&b_3\w\cr
b_4\w&c_3\w&b_2\w\cr
b_3\w&b_2\w&c_4\w\cr\end{matrix}\right]\nnh\nh.
\]%\end{equation}
We have the $\,\mathrm{SO}\hh(3)\,$ freedom of rotating $\,u_2\w,u_3\w,u_4\w$,
leading -- see Lemma~\ref{eqvrt}(b) -- to the same simultaneous rotation of 
the bases (\ref{bas}) of $\,\Lambda\nnh^\pm\nnh$, which allows us to
simultaneously di\-ag\-o\-nal\-ize $\,W^+\nnh$ and $\,W^-\nh$, 
with opposite spectra and, in view of Sect.\,\ref{pp} and~\ref{ft}, proves
Theorem~\ref{first}.

\section{Step four: two double eigen\-val\-ues}\label{td}
We may assume %a special case of (\ref{ord}), namely,
(\ref{lbm}) with $\,\lambda=\mu>0$. By (\ref{set}) and Lemma~\ref{solns}(b),
\[
(W_{\!1212}\w,W_{\!1313}\w,W_{\!1414}\w)
=(W_{\!3434}\w,W_{\!2424}\w,W_{\!2323}\w)=(-\sca/12,\yq-\sca/12,-\yq+\sca/6)
\]
for some $\,\yq$. From (\ref{wij}) -- (\ref{eqv}-a) and (\ref{lbm}) with
$\,\lambda=\mu>0$, for some $\,b_2\w,b_3\w$,
\[
\begin{array}{l}
W_{\!1323}\w=W_{\!1442}\w=b_3\w\hh,\qquad W_{\!1314}\w=W_{\!2342}\w=b_2\w\hh,\\
(W_{\!1213}\w,W_{\!1214}\w,W_{\!1223}\w,W_{\!1224}\w,
W_{\!1334}\w,W_{\!1434}\w,W_{\!2334}\w,W_{\!2434}\w)
=(0,0,0,0,0,0,0,0)\hh.
\end{array}
\]
Thus, 
$\,W\nh$ has in the basis 
$\,u_1\w\wedge u_2\w,\,u_1\w\wedge u_3\w,\,u_1\w\wedge u_4\w,
u_3\w\wedge u_4\w,\,u_4\w\wedge u_2\w,\,u_2\w\wedge u_3\w$ the matrix
\[
\left[\begin{matrix}
-\sca/12&0&0&c_2\w&0&0\cr
0&\yq-\sca/12&b_2\w&0&c_3\w&b_3\w\cr
0&b_2\w&-\yq+\sca/6&0&b_3\w&c_4\w\cr
c_2\w&0&0&-\sca/12&0&0\cr
0&c_3\w&b_3\w&0&\yq-\sca/12&b_2\w\cr
0&b_3\w&c_4\w&0&b_2\w&-\yq+\sca/6\end{matrix}\right]\nnh\nh,
\]
cf.\ (\ref{arb}), and the matrices of
$\,W^\pm\nh$ in the bases (\ref{bas}) of $\,\Lambda\nnh^\pm$ are
\begin{equation}\label{tdb}
\left[\begin{matrix}
\pm c_2\w-\sca/12&0&0\cr
0&\pm c_3\w+\yq-\sca/12&b_2\w\pm b_3\w\cr
0&b_2\w\pm b_3\w&\pm c_4\w-\yq+\sca/6\cr\end{matrix}\right]\nnh\nh.
\end{equation}
We have the $\,\hs\mathrm{SO}\hh(2)\times\mathrm{SO}\hh(2)\,$ 
freedom of independently rotating the pairs $\,u_1\w,u_2\w$ and
$\,u_3\w,u_4\w$. %The final clause of 
Lemma~\ref{eqvrt}(c) now allows us to require that $\,b_2\w=b_3\w=0$, which
proves Theorem~\ref{third}.

Finally, Theorem~\ref{tzero} is a trivial consequence of
Theorems~\ref{first} --~\ref{third} for the multiplicities  
$\,\text{\medbf1111},\,\,\text{\medbf211},\,\,\text{\medbf31},
\,\,\text{\medbf22}\,$ in
(\ref{str}), and of Lemma~\ref{eqvrt}(a) in the case $\,\text{\medbf4}$.

\section{Known geometric examples}\label{kg}
Due to (\ref{iff}), obvious examples of weakly Ein\-stein
Riemannian four-manifolds are provided by those
having $\,\ein=0\,$ (the Ein\-stein ones) or $\,W\nnh\nh=0\,$ and
$\,\sca=0\,$ (the case of con\-for\-mal flatness with zero scalar
curvature). 
There are only three other, quite narrow, classes of examples
known: the EPS space \cite{EPS14}, certain special K\"ah\-ler surfaces 
\cite{derdzinski-euh-kim-park}, and some warped products
\cite[Sect.\,14--15]{derdzinski-park-shin}.

The goal of this section
is to explain where the examples of \cite{derdzinski-euh-kim-park,EPS14} stand
with regard to our algebraic classification. For the examples of
\cite{derdzinski-park-shin}, this question is settled in
\cite[Sect.\,21]{derdzinski-park-shin}.

In the case of the EPS space constructed by 
Euh, Park, and Sekigawa \cite[Example 3.7]{EPS14}, which
-- according to \cite[formula (4.6)]{derdzinski-park-shin} -- 
forms a single lo\-cal-homo\-thety type, $\,R\,$ arises as in
Theorem~\ref{secnd} from any given $\,\lambda>0\,$ and
$\,(\sca,\mu,c_2\w,c_3\w,c_4\w)=(-\nh2\lambda,0,0,0,0)$. Thus,
$\,W^\pm\nnh\nh$ both have the same spectrum
$\,(\lambda/6,\lambda/6,-\nh\lambda/3)$. In fact, using the description
of the curvature tensor $\,R\,$ of the EPS space in 
\cite[Sect.\,3]{derdzinski-euh-kim-park}, and declaring 
our $\,(u_1\w,u_2\w,u_3\w,u_4\w)\,$ and $\,\lambda\,$ %to be 
equal to their $\,(e_1\w,e_3\w,e_4\w,e_2\w)\,$ and $\,2a^2$, one gets 
$\,R_{1212}\w=R_{1313}\w=R_{1414}\w=R_{2323}\w=-R_{2424}\w=-R_{3434}\w
=-\nh\lambda/2$, so that $\,\sca=-\nh2\lambda\,$ and
$\,u_1\w,\dots,u_4\w$ di\-ag\-o\-nal\-ize 
the Ein\-stein tensor $\,\ein\,$ with the eigen\-val\-ues
$\,(-\nnh\lambda,0,0,\lambda)$, while $\,c_2\w=c_3\w=c_4\w=0$.

The construction of weakly Ein\-stein K\"ah\-ler surfaces
in \cite[Sect.\,12]{derdzinski-euh-kim-park} uses orthogonal
vector fields $\,e\nh_i\w$. We choose $\,u\hn_i\w$ related to their 
$\,e\nh_i\w$ by $\,(u_1\w,u_2\w,u_3\w,u_4\w)
=(\hat e_1\w,\hat e_3\w,\hat e_2\w,\hat e_4\w)$, with
$\,\hat e\nh_i\w=e\nh_i\w/|e\nh_i\w|$.
Thus, $\,u_1\w,\dots,u_4\w$ di\-ag\-o\-nal\-ize 
$\,\ein\,$ with some eigen\-val\-ues 
$\,(-\nnh\lambda,-\nnh\lambda,\lambda,\lambda)$, while combining the final
clause of Corollary~\ref{ktype} with formula (12.6)
in \cite[Sect.\,12]{derdzinski-euh-kim-park} we obtain, in our notation,
\[
\begin{array}{l}
(\sca/4,-\yq,\yq-\sca/4)=(c_2\w,c_3\w,c_4\w)
=(W_{\!1234}\w,W_{\!1342}\w,W_{\!1423}\w)\\
\phantom{(\sca/4,-\yq,\yq-\sca/4)=(c_2\w,c_3\w,c_4\w)}
=\,(R_{1234}\w,R_{1342}\w,R_{1423}\w)
=(2p\eta'\nh\theta,-p\eta'\nh\theta,-p\eta'\nh\theta)\hh,
\end{array}
\]
so that $\,\sca\hs=\hs8\hs\xi=8p\eta'\nh\theta$. (The construction uses
positive functions $\,\zeta,\eta\,$ of a real variable and nonzero constants 
$\,p,\theta$.) Formula (12.6) of \cite[Sect.\,12]{derdzinski-euh-kim-park}
also gives, in our notation,
$\,\ric\nh_{11}\w=2[3p\eta'\nh-\zeta\theta\eta'']\hs\theta$. Hence
$\,\ein_{11}\w=\ric\nh_{11}\w-\sca/4
=2[2p\eta'\nh-\zeta\theta\eta'']\hh\theta$. Next, 
$\,\zeta\,$ has the constant derivative $\,-p/\theta$, and the function 
$\,Q=4\hh\zeta\nh\eta\hh\theta^2\nh$,
besides being positive, is only subject to a sec\-ond-or\-der linear
ordinary differential equation 
\cite[formulae (12.4)\,(12.9)]{derdzinski-euh-kim-park}, while
$\,|\hn\lambda|=|\ein_{11}\w\hn|$. Consequently, a suitable choice of
initial data yields the claim made in (\ref{rea}).
\begin{rem}\label{knsqr}We now justify the claim made at the end of the
Introduction. If $\,R\,$ is the Kul\-kar\-ni-No\-mi\-zu square of
$\,\sym$, and $\,\sym\,$ corresponds via $\,\met\,$ to a $\,(1,1)\,$ tensor
$\,B$, then $\,\ric\,$ is analogously related to
$\,(\mathrm{tr}\,B)B-B^2$ and, clearly, $\,\ein=0\,$ when $\,\sym\,$ is a
multiple of $\,\met\,$ or has two mutually
opposite eigen\-values of multiplicities $\,\text{\medbf22}$. In the case 
$\,\text{\medbf31}$, if an orthonormal basis diagonalizes $\,\sym\,$ with
the eigen\-values $\,\delta,\delta,\delta,-\delta$, the corresponding
eigen\-values of $\,\ric\,$ are
$\,\delta^2\nh,\delta^2\nh,\delta^2\nh,-\nh3\delta^2\nh$, and so $\,\sca=0$,
while the only possibly-nonzero components of $\,R\,$ are, with
$\,i,j\in\{1,2,3\}\,$ assumed distinct, those algebraically related to 
$\,R_{ijij}\w=\delta^2$ and $\,R_{i4i4}\w=-\delta^2\nh$, resulting in
$\,W_{\!ijij}\w=R_{ijij}\w-(\delta^2\nh+\delta^2)/2=0\,$ and
$\,W_{\!i4i4}\w=R_{i4i4}\w-(\delta^2\nh-3\delta^2)/2=0$.
\end{rem}

\section*{Acknowledgements}
This work was supported by the National Research Foundation of Korea (NRF)
grant funded by the Korea government (MSIT) (RS-2024-00334956).

%\newpage

\end{document}